\documentclass[12pt]{amsart}%
\usepackage{amscd}
\usepackage{amsmath}
\usepackage[dvips]{graphicx}
\usepackage{amsfonts}
\usepackage{amssymb}%
\theoremstyle{plain}
\newtheorem{theorem}{Theorem}[section]
\newtheorem{corollary}[theorem]{Corollary}
\newtheorem{lemma}[theorem]{Lemma}
\newtheorem{proposition}[theorem]{Proposition}
\theoremstyle{definition}

\newtheorem{definition}[theorem]{Definition}
\theoremstyle{remark}

\newtheorem*{claim}{Claim}

\newtheorem{remark}[theorem]{Remark}

\renewcommand{\emptyset}{\varnothing}
\begin{document}
\title[von Neumann algebra of Baumslag group embeds into $R^{\omega}$]{The von Neumann algebra of the non-residually finite Baumslag group
$\left\langle a,b\mid ab^{3}a^{-1}=b^{2}\right\rangle$ embeds into $R^{\omega}$}
\author{Florin R\u{a}dulescu}
\address{Department of Mathematics\\
The University of Iowa\\
Iowa City, Iowa 52242, U.S.A.}

\begin{abstract}
In this paper we analyze the structure of some sets of non-commutative moments
of elements in a finite von Neumann algebra $M$. If the fundamental group of
$M$ is $\mathbb{R}_{+}\backslash\{0\}$, then the moment sets are convex, and
if $M$ is isomorphic to $M\otimes M$, then the sets are closed under pointwise
multiplication. We introduce a class of discrete groups that we call
hyperlinear. These are the discrete subgroups (with infinite conjugacy
classes) of the unitary group of $R^{\omega}$. We prove that this class is
strictly larger than the class of (i.c.c.) residually finite groups; it
contains the Baumslag group $\left\langle a,b\mid ab^{3}a^{-1}=b^{2}%
\right\rangle $. This leads to a previously unknown (non-hyperfinite) type
$\mathrm{II}_{1}$ factor that can be embedded in $R^{\omega}$. This is
positive evidence for Connes's conjecture that any separable type
$\mathrm{II}_{1}$ factor can be embedded into $R^{\omega}$.

\end{abstract}
\maketitle

\section{Introduction}

In this paper we analyze the structure of sets of (non-commutative) moments
$\tau\left(  x_{1}\cdots x_{n}\right)  $ of variables $x_{1},x_{2},\dots
,x_{n}$ in a type $\mathrm{II}_{1}$ factor $M$. We analyze the structure of
these sets, for the case of projections and unitaries (see also \cite{Ra} for
odd moments of selfadjoint elements). While the understanding of these
structures is far from being complete, we prove that any discrete (i.c.c.)
group $\Gamma$ that can be faithfully embedded into the unitary group of
$R^{\omega}$ has the property that $\mathcal{L}(\Gamma)\subseteq R^{\omega}$.

By using techniques pertaining to free probability we prove that the Baumslag
group $\Gamma=\left\langle a,b\mid ab^{3}a^{-1}=b^{2}\right\rangle $, which is
non-residually finite, embeds (faithfully) into $\mathcal{U}(R^{\omega})$.
Note also that by \cite{Che}, the algebra $\mathcal{L}(\Gamma)$ does not have
property $\Gamma$. In particular, it is non-hyperfinite. (We are indebted to
P. de la Harpe for bringing this to our attention.)

This is positive evidence towards Connes's conjecture that any separable
$\mathrm{II}_{1}$ factor is embedded into $R^{\omega}$.

In the second part of the paper we analyze the structure of the sets of first
and second order of moments $\tau\left(  e_{i}\right)  $, $\tau\left(
e_{i}e_{j}\right)  $ of finite sets of projections $e_{1},\dots,e_{n}$ in a
type $\mathrm{II}_{1}$ factor. By work of Kirchberg \cite{Ki}, if the closure
of the sets of first and second order of moments of unitaries is independent
of the type $\mathrm{II}_{1}$ factor considered, then Connes's conjecture
should be true.

It is obvious how to translate this statement in terms of projections. By
using methods from \cite{Ra}, it follows that the corresponding set of
moments, i.e., the set of moments (of order $1$ and $2$) of projections, is
convex and multiplicative for a $\mathrm{II}_{1}$ factor $M$ such that $M\cong
M\otimes M$ and $\mathcal{F}(M)=\mathbb{R}_{+}\backslash\{0\}$. We will also
analyze the structure of faces of these sets, which gives some additional data
on the geometric structure of these sets.

This work was supported by NSF grant DMS99-70486.\medskip

\emph{Definitions and Notations}: We recall that for a von Neumann algebra
$M$, the unitary group is denoted by $\mathcal{U}(M)$, while $\mathcal{P}(M)$
stands for the set of selfadjoint projections.

For a type $\mathrm{II}_{1}$ factor $M$, the fundamental group $\mathcal{F}%
(M)$ of $M$ \cite{MuN} \cite{vN} is defined as the multiplicative group of all
$t>0$, such that $M_{t}\cong M$.

For $\Gamma$ a countable discrete group, with infinite conjugacy classes
(briefly i.c.c.), the algebra $\mathcal{L}(\Gamma)$ is the weak closure of the
group algebra $\mathbb{C}(\Gamma)$ embedded (via left regular representation)
in $B(\ell^{2}(\Gamma))$.

If $\omega$ is an ultrafilter on $\mathbb{N}$, then following \cite{McD} and
\cite{Co}, one defines for any $\mathrm{II}_{1}$ factor $M$ the ultrafilter
product $M^{\omega}$, obtained via G.N.S. construction, by defining the trace
of an element $(x_{n})$ in the infinite product of copies of $M$ to be
$\lim_{n\rightarrow\omega}$ $\tau(x_{n})$ (in the hypothesis that
$\sup\left\|  x_{n}\right\|  <\infty$). We refer to Connes's \cite{Co} paper
on injectivity for full details on this construction.

Finally, we recall a construction from \cite{Po} (see also \cite{Vo}).
Consider two von Neumann algebras $N_{1}$, $N_{2}$ which have a common
subalgebra $B$, containing the unit. Also assume that the algebras $N_{i}$
have faithful traces whose restriction to $B$ coincides. We assume that we are
given conditional expectations $E_{i}$ from $N_{i}$ onto $B$ that are trace preserving.

The trace on the reduced amalgamated free product von Neumann algebra
$C_{1}\ast_{B}C_{2}$ is defined by the requirement that a product
\linebreak$c_{1,1}c_{2,1}c_{1,2}c_{2,2}c_{1,3}c_{2,3}\cdots$, where $c_{1,i}$
belongs to $C_{1}$ and $c_{2,i}$ belongs to $C_{2}$, has zero trace if
$\operatorname*{Id}-E_{B}(c_{ij})\neq0$ for all $i=1,2$, $j=1,2,\dots$.

\section{\label{Mom}Moments of unitaries}

In this section we define some sets of non-commutative moments of unitaries
$\tau\left(  u_{1}u_{2}\cdots u_{p}\right)  $ associated with a type
$\mathrm{II}_{1}$ factor $M$. We will use these sets to check that for any
discrete i.c.c.\ group $\Gamma$ that can be embedded into the unitary group of
$R^{\omega}$, then also $\mathcal{L}(\Gamma)$ can be embedded into $R^{\omega
}$ (as a unital subfactor).

First we consider the set of all possible embeddings (up to order $N$) of a
group-like algebra.

By $\mathcal{V}_{n,p}$ we denote the set of all indices $(i_{1},i_{2}%
,\dots,i_{k})$ with $1\leq k\leq p$ and $i_{1},i_{2},\dots,i_{k}$ in
$\{1,2,\dots,n\}$. By $u_{I}$ we denote the product $u_{i_{1}}u_{i_{2}}\cdots
u_{i_{k}}$ if $I=(i_{1},i_{2},\dots,i_{k})$.

\begin{definition}
\label{definition2.1}Let $M$ be a separable type $\mathrm{II}_{1}$ factor and
consider the following subset $($of $\{0,1\}^{2^{n^{p}}})$, denoted by
$K_{M}^{n,p}$. We define $K_{M}^{n,p}$ by requiring that $(\varepsilon
_{I})_{\left|  I\right|  \leq p}$ belongs to $K_{M}^{n,p}$ if and only if
there exist unitaries $u_{1},u_{2},\dots,u_{n}$ in $\mathcal{U}(M)$ such that
$\tau(u_{I})=1$ or $0$, and $\varepsilon_{I}=\tau(u_{I})$, for $\left|
I\right|  \leq p$. Here $I$ is an index set $(i_{1},i_{2},\dots,i_{k})$, with
$ij\in\{1,2,\dots,n,\}$, $1\leq k\leq p$ and $\left|  I\right|  =k$.
\end{definition}

\begin{remark}
\label{RemMomNew.2}If $M=\mathcal{L}\left(  \Gamma\right)  $ and $\Gamma$ has
non-solvable world problem, then there are $2^{n^{p}}$-uples of $0$'s and
$1$'s about which it might be undecidable whether they belong to $K_{M}^{n,p}$.
\end{remark}

The above definition might be restrictive for some purposes because it
requires that $\tau(u_{I})$ is either $0$ or $1$. In fact $\tau(u_{I})=1$ is
equivalent to the fact that $u_{I}=1$. (Here by $u_{I}$ we mean the product
$u_{i_{1}}u_{i_{2}}\cdots u_{i_{k}}$ if $I=(i_{1},i_{2},\dots,i_{k})$.)

\begin{definition}
\label{definition2.2}Let $M$ be a finite von Neumann algebra and let for fixed
integers $n,p$,
\begin{multline*}
L_{M}^{n,p}=\{(\tau(u_{I})=\tau(u_{i_{1}}u_{i_{2}}\cdots u_{k}))_{I\in
\mathcal{V}_{n,p}}\mid\text{for all}\\
\text{ unitaries }u_{1},u_{2},\dots,u_{n}\text{ in }\mathcal{U}(M)\}.
\end{multline*}

\end{definition}

The following properties are easy to observe \textup{(}see also \cite{Ra}%
\textup{).}

\begin{proposition}
\label{proposition2.3}

\begin{enumerate}
\item[$\mathrm{(a)}$] Let $M_{1}$, $M_{2}$ be finite von Neumann algebras.
Denote by $\odot$, the pointwise product on $\mathbb{R}^{s}$, for all $s$.
Then for all positive integers $n,p$,
\[
L_{M_{1}}^{n,p}\odot L_{M_{2}}^{n,p}\subseteq L_{M_{1}\otimes M_{2}}%
^{n,p},\qquad K_{M_{1}}^{n,p}\odot K_{M_{2}}^{n,p}\subseteq K_{M_{1}\otimes
M_{2}}^{n,p}.
\]

\item[$\mathrm{(b)}$] In particular if $M$ is such that $M\cong M\otimes M$,
then $K_{M}^{n,p}$ and $L_{M}^{n,p}$ are closed under pointwise multiplication.

\item[$\mathrm{(c)}$] If $\lambda\in\mathcal{F}(M)$, then $\lambda L_{M}%
^{n,p}+(1-\lambda)L_{M}^{n,p}\subseteq L_{M}^{n,p}$ for all integers
$n,p\geq1$. In particular if $\mathcal{F}(M)=\mathbb{R}_{+}\backslash\{0\}$,
then $L_{M}^{n,p}$ is convex.

\item[$\mathrm{(d)}$] $\overline{L_{M}^{n,p}}\subseteq L_{M^{\omega}}^{n,p}$
and $L_{M^{\omega}}^{n,p}$ is closed in the product topology of $\mathbb{R}%
^{\left|  \mathcal{V}_{n,p}\right|  }$.

\item[$\mathrm{(e)}$] In particular if $M\cong M\otimes M$ and $(\lambda
_{I})_{I\in\mathcal{V}_{n,p}}$ is an element in $L_{M}^{n,p}$ such that either
$\left|  \lambda_{I}\right|  <1$ or $\lambda_{I}=1$, then by replacing the
components in $L_{M}^{n,p}$ which are not $1$, by zero, we obtain an element
in $K_{M^{\omega}}^{n,p}$.

\item[$\mathrm{(f)}$] Let $\Phi_{u_{1}}$ be the operation on $L_{M}^{n,p}$
which replaces in $(\lambda_{I})_{I\in L_{M}^{n,p}}$ any monomial $\lambda
_{I}=\tau(u_{I})$, corresponding to a nonzero total power of $u_{1}$, by zero.
Assume $\mathcal{F}(M)=\mathbb{R}_{+}\backslash\{0\}$. Then
\[
\Phi_{u_{1}}\left(  L_{M}^{n,p}\right)  \subseteq L_{M}^{n,p}.
\]

\item[$\mathrm{(g)}$] For all separable $\mathrm{II}_{1}$ factors $M$ we have
that $L_{M}^{n,p}\supseteq L_{R}^{n,p}$. Moreover if $M\subseteq R^{\omega}$,
then $\overline{L_{M}^{n,p}}=\overline{L_{R}^{n,p}}$ \textup{(}and similarly
for $K$\textup{).}
\end{enumerate}
\end{proposition}

Open Question: Does $K_{\mathcal{L}(\mathrm{SL}_{3}(\mathbb{Z}))}%
^{n,p}\subseteq K_{\mathcal{L}(F_{\infty})}^{n,p}$ for all $n,p$?

The proof of properties $\mathrm{(a)}$--$\mathrm{(c)}$ is obvious and is
basically contained in \cite{Ra}. To check property $\mathrm{(d)}$ we have
only to verify that $L_{M^{\omega}}^{n,p}$ is closed. But if $(\lambda
_{I})_{I\in\mathcal{V}_{n,p}}$ is an accumulation point for $L_{M}^{n,p}$,
then take unitaries $(u_{i}^{s})_{i=1}^{n}$, for all $s$, such that
$\lim_{s\rightarrow\infty}\tau(u_{I}^{s})=\lambda_{I}$, $I\in\mathcal{V}%
_{n,p}$. Then $u_{i}=(u_{i}^{s})_{s\in\mathbb{N}}$ are unitaries in
$R^{\omega}$, whose non-commutative moments give $(\lambda_{I})_{I\in
\mathcal{V}_{n,p}}$.

Property $\mathrm{(e)}$ follows from properties $\mathrm{(a)}$ and
$\mathrm{(c)}$. Property $\mathrm{(f)}$ follows by convexity, and integration
over $\theta$, where gauging $u_{1}$ by $e^{2\pi i\theta}$. Property
$\mathrm{(g)}$ is obvious.

\begin{proposition}
\label{proposition2.4}Let $\Gamma$ be a discrete i.c.c.\ group that embeds
faithfully into the unitary group of $R^{\omega}$. Then $\mathcal{L}%
(\Gamma)\subseteq R^{\omega}$.
\end{proposition}

\begin{proof}
Fix $n,p$ and let $u_{1},u_{2},\dots,u_{n},\dots$ be a system of generators of
$\Gamma$. Let $\varepsilon_{I}$ be the traces of $(u_{I})_{I\in\mathcal{V}%
_{n,p}}$ in the left regular representation of $\mathcal{L}(\Gamma)$. By
hypothesis there exist unitaries $v_{1},v_{2},\dots,v_{n}$ in $R^{\omega}$
such that $\left|  \tau(v_{I})\right|  <1$ if $u_{I}\neq1$ in $\Gamma$ and
$v_{I}=1$ if $u_{I}=1$ in $\Gamma$. Let $(\alpha_{I})^{s}=\tau(v_{I}^{\otimes
s})=\tau(v_{I})^{s}$. Then $(\alpha_{I})_{I\in\mathcal{V}_{n,p}}^{s}$ belongs
to $\overline{L_{R}^{n,p}}$ and hence so does the limit
\[
\varepsilon_{I}=\underset{s\rightarrow\infty}{\lim}(\alpha_{I})^{s}%
\text{\qquad for }I\in\mathcal{V}_{n,p}.
\]
Thus $(\varepsilon_{I})_{I\in\mathcal{V}_{n,p}}$ $\in\overline{K_{R}^{n,p}}$
for all $n,p$. Hence $\mathcal{L}(\Gamma)\subseteq R^{\omega}$.
\end{proof}

\begin{definition}
\label{definition2.5}We call an i.c.c.\ group $\Gamma$ hyperlinear if $\Gamma$
embeds faithfully into $\mathcal{U}(R^{\omega})$. Clearly any residually
finite group is hyperlinear. The class of hyperlinear groups is obviously
closed under free products.
\end{definition}

\begin{theorem}
\label{theorem2.6}The class of hyperlinear groups is strictly larger than the
class of residually finite groups. More precisely, the Baumslag group
$\left\langle a,b\mid ab^{3}a^{-1}=b^{2}\right\rangle $ \cite{Ba} \cite{Ma} is
hyperlinear and non-residually finite. \textup{(}Note that by \cite{Che},
$\mathcal{L}(\Gamma)$ does not have property $\Gamma$.\textup{)}
\end{theorem}

\begin{proof}
We divide the proof into several steps. We construct first an approximate
embedding of the relation $ab^{3}a^{-1}=b^{2}$ into $M_{n}(\mathbb{C})$. We
then take the free amalgamated product (over $b^{2}$) by a unitary that
commutes with $b^{2}$ and perturb $a$ with this unitary. This gives an
approximate embedding of $\Gamma$ into the nonscalar unitaries in some free
product algebras. Since these algebras are embeddable into $R^{\omega}$
\cite{Vo} \cite{Wa}, the result will follow, by Proposition
\ref{proposition2.4}.\medskip

\emph{Step} I. Construction of an approximate embedding.

There exist unitaries $v_{n}$, $b_{n}$ (of zero trace) in $M_{6n}(\mathbb{C})$
with the following properties:\medskip

Property 1)%
\[
\left\|  v_{n}b_{n}^{3}v_{n}-b_{n}^{2}\right\|  _{\infty}\leq\frac{K_{1}}%
{n}\text{for a universal constant }K_{1}.
\]

Property 2)%

\begin{quotation}%
Denote by $B_{0}^{n}$ the abelian algebra generated by $b_{n}^{2}$ and let
$E_{B_{0}^{n}}$ be the corresponding conditional expectation. Let $\Phi
^{n}=\operatorname*{Id}-E_{B_{0}^{n}}$. Then for all $\alpha\in\{\pm1,\pm
2\}$,
\[
\left\|  \Phi(v_{n}b_{n}^{\alpha}v_{n}^{\ast})\right\|  _{2}\geq K_{2}%
\]
for a universal constant $K_{2}$. Here $\left\|  \;\;\right\|  _{2}$ is the
normalized Hilbert-Schmidt trace on matrices.%
\end{quotation}%

Property 3)%
\[
E_{B_{0}^{n}}(b_{n}^{\pm1})=0,\;E_{B_{0}^{n}}\left(  vb^{\alpha}\right)
=0,\;E_{B_{0}^{n}}\left(  b^{\alpha}v\right)  =0,\qquad\alpha\in
\mathbb{N},\/\alpha\neq0.
\]

We describe first the construction of the unitaries $v_{n}$, $b_{n}$.

Let $e_{0},e_{1},\dots,e_{6n-1}$ be the diagonal algebra of $M_{6n-1}%
(\mathbb{C})$, and for convenience we think of $e_{k}$ as being identified
with $\chi_{\lbrack\frac{k}{6n},\frac{k+1}{6n})}^{{}}$, $k=0,1,\dots,6n-1$.

Let $f_{k}$, for $k=0,1,2,\dots,n-1$, be the projection $\chi_{\lbrack
\frac{3k}{6n},\frac{3k+3}{6n})+\{0,\frac{1}{2}\}}^{{}}$ and let
\[
g_{k}=\chi_{\lbrack\frac{2k}{6n},\frac{2k+2}{6n})+\left\{  0,\frac{1}{3}%
,\frac{2}{3}\right\}  }^{{}}.
\]

Let $b_{n}$ be the unitary defined by
\[
b_{n}=\sum_{k=0}^{6n-1}e^{\frac{2\pi ik}{6n}}e_{k}.
\]
Let $v_{n}$ be a unitary such that $v_{n}^{\ast}f_{k}=g_{k}v_{n}^{\ast}$,
$k=0,1,2,\dots,n-1$, and such that
\begin{align*}
(\operatorname*{Ad}v_{n}^{\ast})(e_{3k+\varepsilon}) &  =e_{2k+\varepsilon},\\
(\operatorname*{Ad}v_{n}^{\ast})(e_{3k+3n+\varepsilon}) &
=(e_{2k+4n+\varepsilon})
\end{align*}
for all $k=0,1,2,\dots,n-1$, $\varepsilon=0,1$. (See Fig.\ \ref{Figvnstar}.)

\begin{figure}[ptb]
\setlength{\unitlength}{360bp} \begin{picture}(1,1)
\put(0,0){\includegraphics[bb=0 0 360 360,height=360bp,width=360bp]{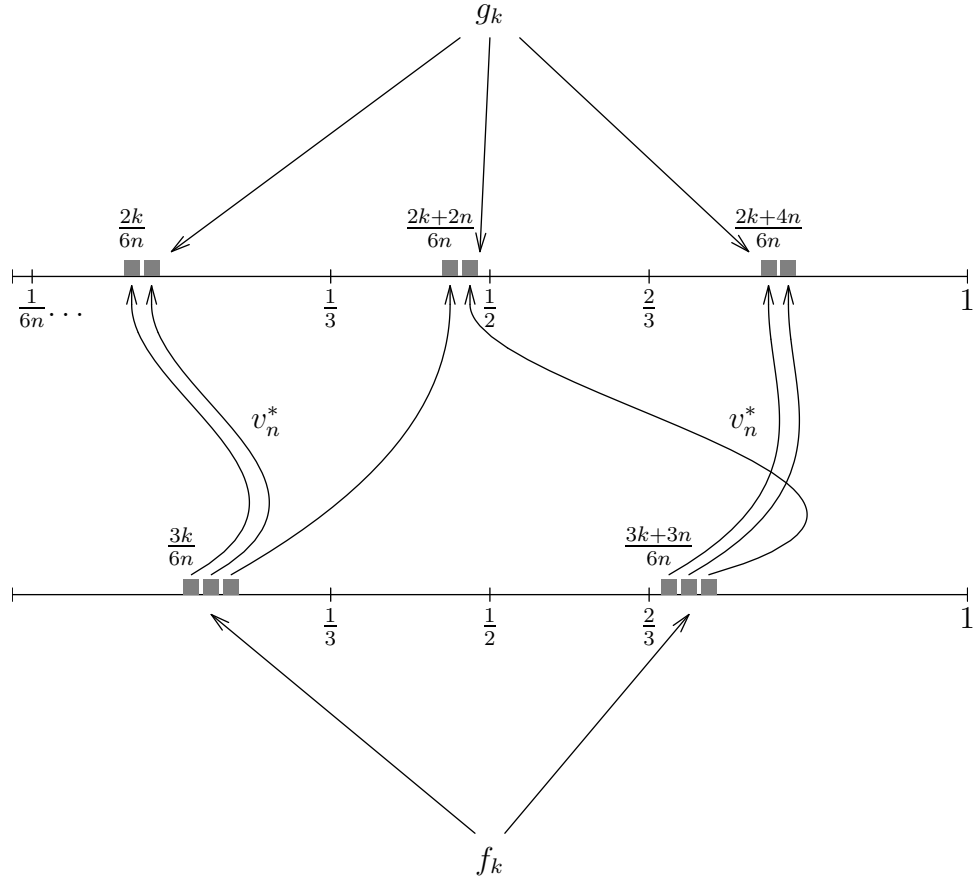}}
\put(0.5,0.0694444){\makebox(0,0)[t]{$f_{k}$}}
\put(0.5,0.9305556){\makebox(0,0)[b]{$g_{k}$}}
\put(0.0208333,0.6527778){\makebox(0,0)[t]{$\frac{1}{6n}$\rlap{$\dots$}}}
\put(0.3333333,0.6527778){\makebox(0,0)[t]{$\frac{1}{3}$}}
\put(0.5,0.6527778){\makebox(0,0)[t]{$\frac{1}{2}$}}
\put(0.6666667,0.6527778){\makebox(0,0)[t]{$\frac{2}{3}$}}
\put(1,0.6527778){\makebox(0,0)[t]{$1$}}
\put(0.3333333,0.3194444){\makebox(0,0)[t]{$\frac{1}{3}$}}
\put(0.5,0.3194444){\makebox(0,0)[t]{$\frac{1}{2}$}}
\put(0.6666667,0.3194444){\makebox(0,0)[t]{$\frac{2}{3}$}}
\put(1,0.3194444){\makebox(0,0)[t]{$1$}}
\put(0.125,0.697917){\makebox(0,0)[b]{$\frac{2k}{6n}$}}
\put(0.447917,0.697917){\makebox(0,0)[b]{$\frac{2k+2n}{6n}$}}
\put(0.791667,0.697917){\makebox(0,0)[b]{$\frac{2k+4n}{6n}$}}
\put(0.177083,0.364583){\makebox(0,0)[b]{$\frac{3k}{6n}$}}
\put(0.677083,0.364583){\makebox(0,0)[b]{$\frac{3k+3n}{6n}$}}
\put(0.25,0.5){\makebox(0,0)[bl]{$v_{n}^{*}$}}
\put(0.75,0.5){\makebox(0,0)[bl]{$v_{n}^{*}$}}
\end{picture}
\caption{Description of $v_{n}^{\ast}$}%
\label{Figvnstar}%
\end{figure}

For $\varepsilon=2$, we define $v_{n}^{\ast}$ by the requirement
\begin{align*}
\operatorname*{Ad}v_{n}^{\ast}(e_{3k+2}) &  =e_{2k+2n},\\
\operatorname*{Ad}v_{n}^{\ast}(e_{3k+3n+2}) &  =e_{2k+2n+1}.
\end{align*}
Observe that the definition of $b_{n}$, $v_{n}$ implies, for a universal
constant $K_{1}$, that%
\begin{equation}
\left\Vert b_{n}^{2}f_{k}-e^{\frac{2\pi ik}{n}}f_{k}\right\Vert _{\infty}%
\leq\frac{K_{1}}{n},\label{eq1}%
\end{equation}%
\begin{equation}
\left\Vert b_{n}^{3}g_{k}-e^{\frac{2\pi ik}{n}}g_{k}\right\Vert _{\infty}%
\leq\frac{K_{1}}{n},\qquad k=0,1,\dots,n-1.\label{eq2}%
\end{equation}
Moreover, $v_{n}^{\ast}f_{k}^{{}}=g_{k}^{{}}v_{n}^{\ast}$, and%
\begin{equation}
(v_{n}b_{n}v_{n}^{\ast})\cdot e_{3k+\varepsilon+\alpha\cdot3n}=e^{2\pi
i(\frac{2k+\varepsilon}{6n}+\alpha\cdot\frac{2}{3})}e_{3k+\varepsilon
+\alpha\cdot3n}\text{\qquad for }\varepsilon=0,1,\;\alpha=0,1.\label{eq3}%
\end{equation}
In the remaining case, we have that
\begin{equation}
(v_{n}b_{n}v_{n}^{\ast})\cdot e_{3k+2+\alpha\cdot3n}=e^{2\pi i(\frac
{2k+\alpha}{6n}+\frac{1}{3})}e_{3k+2+\alpha\cdot3n}\text{\qquad for }%
\alpha=0,1.\label{eq4}%
\end{equation}

We now proceed to the proof of the properties 1), 2), 3).

Since $v_{k}^{\ast}f_{k}^{{}}=g_{k}^{{}}v_{k}^{\ast}$, $v_{k}g_{k}=f_{k}v_{k}%
$, $v_{n}^{{}}g_{k}^{{}}v_{n}^{\ast}=f_{k}^{{}}$, and $f_{k},g_{k}$ commute
with $b_{n}$, it follows that
\begin{align*}
\left\|  v_{n}b_{n}^{3}v_{n}^{\ast}-b_{n}^{2}\right\|  _{\infty}  &
=\underset{k=0,1\dots,n-1}{\max}\left\|  (v_{n}b_{n}^{3}v_{n}^{\ast}-b_{n}%
^{2})f_{k}\right\| \\
&  =\underset{k=0,1\dots,n-1}{\max}\left\|  v_{n}(b_{n}^{3}g_{k}-e^{\frac{2\pi
ik}{n}}g_{k})v_{n}^{\ast}-(b_{n}^{2}f_{k}-e^{\frac{2\pi ik}{n}}f_{k})\right\|
_{\infty}.
\end{align*}
But this quantity is less than $\frac{K_{1}}{n}$ by (\ref{eq1}), (\ref{eq2}).
This completes the proof of Property 1.

To prove Property 2 we need to describe first $E_{B_{0}^{n}}$. But it is
obvious that
\[
E_{B_{0}^{n}}\left(  \sum_{k=0}^{6n-1}\lambda_{k}e_{k}\right)  =\sum
_{k=0}^{3n-1}\frac{1}{2}\left(  \lambda_{k}+\lambda_{k+3n}\right)
(e_{k}+e_{k+3n}).
\]

Consequently
\[
\Phi_{n}\left(  \sum_{k=0}^{6n-1}\lambda_{k}e_{k}\right)  =\sum_{k=0}%
^{3n-1}\frac{\lambda_{k}-\lambda_{k+3n}}{2}e_{k}+\sum_{k=0}^{3n-1}%
\frac{\lambda_{k+3n}-\lambda_{k}}{2}e_{k+3n}.
\]
We use the above formula for $v_{n}^{{}}b_{n}^{{}}v_{n}^{\ast}$ and use
(\ref{eq3}) and (\ref{eq4}). We take $P_{n}=\sum_{k=0}^{n-1}\chi_{\lbrack
\frac{3k}{6n},\frac{3k+1}{6n})}^{{}}=\sum_{k=0}^{n-1}e_{3k}$. Then $P_{n}$ has
trace $\frac{1}{6}$. Since $\left(  v_{n}^{{}}b_{n}^{{}}v_{n}^{\ast}\right)
e_{3k}=e^{2\pi i\frac{k}{3n}}e_{3k}$ and $\left(  v_{n}^{{}}b_{n}^{{}}%
v_{n}^{\ast}\right)  e_{3k+3n}=e^{2\pi i\left(  \frac{k}{3n}+\frac{2}%
{3}\right)  }e_{3k+3n}$ the above formula shows that
\[
P_{n}\Phi_{n}(vbv^{\ast})=P_{n}(vbv^{\ast}-E_{B_{0}^{n}}(vbv^{\ast}%
))=\sum_{k=0}^{n-1}e^{2\pi i\frac{k}{3n}}\left(  1-e^{2\pi i\frac{2}{3}%
}\right)  e_{3k}.
\]
Hence
\[
\left\Vert P_{n}\Phi_{n}(vbv^{\ast})\right\Vert _{2}^{2}>\left\vert 1-e^{2\pi
i\frac{2}{3}}\right\vert \left\Vert P_{n}\right\Vert _{2}^{2}=\frac{1}%
{6}\left\vert 1-e^{2\pi i\frac{2}{3}}\right\vert .
\]
The computations for $vb_{n}^{\pm2}v^{\ast}$, $vb_{n}^{-1}v^{\ast}$ are
similar, eventually the factor $\frac{2}{3}$ being replaced by $\frac{4}{3}$
or $-\frac{2}{3}$. This completes the proof of Property 2.

It is obvious that $E_{B_{0}^{n}}(vb_{n}^{\alpha})$, $E_{B_{0}^{n}}%
(b_{n}^{\alpha}v)$, $\alpha\neq0$, and $E_{B_{0}^{n}}(b_{n}^{\pm1})$ are
vanishing.\medskip

\emph{Step} II. In this step we consider the amalgamated free product of the
algebra $\{v_{n},b_{n}\}^{\prime\prime}$ described above and $\mathcal{L}%
(\mathbb{Z})\otimes B_{0}^{n}$. The amalgamated free product is considered
over $B_{0}^{n}$ (the von Neumann algebra generated by $b_{n}^{2}$).

Let $\mathcal{L}(\mathbb{Z})$ have the canonical generator $a_{1}$, a Haar
unitary. $\mathcal{L}(\mathbb{Z})$ is endowed with the standard trace.
Consider the algebra
\[
\mathcal{A}_{n}=(\mathcal{L}(\mathbb{Z})\otimes B_{0}^{n})\ast_{B_{0}^{n}%
}\{v_{n},b_{n}\}^{\prime\prime}%
\]
with the canonical amalgamated free product trace (see section on definitions,
\cite{Po}, and \cite{Vo}).

By \cite{Ra} (see also \cite{Dy}, \cite{Shly}), we have that $\mathcal{A}_{n}$
is a free group factor. Using the ultrafilter construction (\cite{Co},
\cite{McD}), we construct algebras $\mathcal{A}^{\omega}$, $B_{0}^{\omega}$
and $B^{\omega}$ consisting of bounded sequences of elements in the algebras
$\mathcal{A}_{n}$, $B_{0}^{n}$ and $\left\{  b_{n}\right\}  ^{\prime\prime}$
respectively. It is obvious that for $x=\left(  x_{n}\right)  _{n}$ in
$\mathcal{A}^{\omega}$ we have
\[
E_{B_{0}^{\omega}}(\left(  x_{n}\right)  _{n})=\left(  E_{B_{0}^{n}}%
(x_{n})\right)  _{n}.
\]
Let $b$ be the unitary element $\left(  b_{n}\right)  _{n}\in B^{\omega}$ and
let $B_{0}$ be the (abelian) von Neumann algebra generated by $b$. Let $v$ be
the unitary $v=\left(  v_{n}\right)  _{n}\in\mathcal{A}^{\omega}$. We identify
$\mathcal{L}\left(  \mathbb{Z}\right)  \subseteq\mathcal{A}^{\omega}$ with
constant sequences with elements in $\mathcal{L}\left(  \mathbb{Z}\right)  $.

Since, by \cite{Wa}, \cite{Vo}, any type $\mathrm{II}_{1}$ free group factor
embeds into $R^{\omega}$, we obtain that the algebra $\mathcal{A}^{\omega}$ is
embedded into $R^{\omega}$ and hence
\[
\mathcal{A}=\mathcal{L}(\mathbb{Z})\otimes B_{0}\ast_{B_{0}}\{v,b\}^{\prime
\prime}\subseteq R^{\omega}.
\]
The trace on $\mathcal{A}$ is the amalgamated free product trace and coincides
with the restriction of the trace on $R^{\omega}$. Let $\Phi$ be the identity
minus the conditional expectation $E_{B_{0}}$ from $\mathcal{A}$ (or
$\mathcal{A}^{\omega}$) onto $B_{0}$. The following properties hold true:

\begin{enumerate}
\item[1)] $v,b$ are Haar unitaries, $vb^{3}v^{\ast}=b^{2}$.

\item[2)] \raggedright$\left\Vert \Phi\left(  vb^{\alpha}v^{\ast}\right)
\right\Vert _{2}\geq\frac{1}{6}$, $\alpha\in\{\pm1,\pm2\}$; $E_{B_{0}}%
(vb^{k})=0$, $k\in\mathbb{N}$, $E_{B_{0}}(b^{\pm1})=0$.
\end{enumerate}

To prove Property 2, note that $\left\Vert E_{B_{0}^{\omega}}\left(
vb^{\alpha}v^{\ast}\right)  \right\Vert _{2}\leq\frac{1}{6}$ because of the
corresponding property for $v_{n}b_{n}^{\alpha}v_{n}^{\ast}$. Since
$B_{0}\subseteq B_{0}^{\omega}$ it follows also that $\left\Vert E_{B_{0}%
}\left(  vb^{\alpha}v^{\ast}\right)  \right\Vert _{2}\leq\frac{1}{6}$.\medskip

Let $a_{1}$ be the standard generator of $\mathcal{L}(\mathbb{Z})$ and let
$A=a_{1}v$, $B=b$.\medskip

\emph{Step} III. Let $A,B$ be the unitaries defined in Step II. Clearly
\linebreak$AB^{3}A^{-1}=B^{2}$, as $a_{1}$ commutes with $b^{2}$. Let
\[
W=A^{\alpha_{1}}B^{\beta_{1}}A^{\alpha_{2}}\cdots A^{\alpha_{n}}B^{\beta
_{n+1}}%
\]
be a word in $A,B$ such that $\beta_{1}\neq0,\dots,\beta_{n}\neq0$,
$\alpha_{2}\neq0,\dots,\alpha_{n}\neq0$. Consider the following assumption on
the sequence of the indices $\alpha_{i}$.\medskip

\noindent\emph{Assumption} on $W=A^{\alpha_{1}}B^{\beta_{1}}A^{\alpha_{2}%
}\cdots A^{\alpha_{n}}B^{\beta_{n+1}}$. One of the following possibilities
occurs (about consecutive indices):

\begin{enumerate}
\item[A1)] Either $\alpha_{i}$, $\alpha_{i+1}$ are both positive or negative
(except for the case when $\alpha_{1}=0)$.

\item[A2)] If $\alpha_{i}<0$, $\alpha_{i+1}>0$, then $\beta_{i}\in\{\pm1\}$.

\item[A3)] If $\alpha_{i}>0$, $\alpha_{i+1}<0$, then $\beta_{i}\in\{\pm
1,\pm2\}$.
\end{enumerate}

\begin{claim}
\label{claim2}If the word $W$ is subject to the conditions A1, A2, A3, then
$W$ is not a multiple of a scalar (and hence $\left|  \tau(W)\right|  <1$).
\end{claim}

\begin{proof}
[Proof of the claim in Step \textup{III}]We use the following property of an
\linebreak amalgamated free product $\mathcal{B}=E\ast_{C}F$ where $E$, $F$ are
finite algebras with faithful traces $\tau_{1}$, $\tau_{2}$ whose restrictions
coincide on the common unital subalgebra $C$.\renewcommand{\qed}{}
\end{proof}

Assume $w=e_{1}f_{1}e_{2}f_{2}\cdots e_{n}f_{n+1}$ is a word in $E\ast_{C}F$,
$e_{i}\in E$, $f_{i}\in F$, such that $\operatorname*{Id}-E_{C}(f_{1})$,
$\operatorname*{Id}-E_{C}(e_{2})\neq0$, $\dots$, $\operatorname*{Id}%
-E_{C}(e_{n})\neq0$. Then $w$ is not a scalar multiple of the identity. This
follows for example from the construction in \cite{Po}.

Then for the word $W=A^{\alpha_{1}}B^{\beta_{1}}A^{\alpha_{2}}\cdots
A_{g}^{\alpha_{n+1}}$, we use the fact that $A=a_{1}v$, $B=b$. Since
\[
E_{B_{0}}(b^{\theta}v^{\ast}),\text{ }E_{B_{0}}(vb^{\theta})
\]
are always zero for all $\theta$, the only instances in the product in $W$
where we could have elements with $\operatorname*{Id}-E_{B_{0}}$ nonzero are
in subsequences of the form
\[
\cdots avb^{\pm\alpha}v^{\ast}a\cdots,\qquad\alpha\in\{1,2\}\text{ (in }\cdots
AB^{\pm\alpha}A^{-1}\cdots\text{),}%
\]
or
\[
\cdots v^{\ast}a^{\ast}b^{\pm1}av\cdots\text{\qquad(in }\cdots A^{-1}B^{\pm
1}A\cdots\text{).}%
\]
But in these cases $\Phi=\operatorname*{Id}-E_{B_{0}}$ applied to the elements
$vb^{\pm1}v^{\ast}$, $vb^{\pm2}v^{\ast}$, and $b^{\pm1}$ is nonzero. The
remaining two cases correspond to subsequences involving $A^{n}B^{\theta}%
A^{m}$ with $\theta\neq0$ and $n$, $m$ both strictly positive or both strictly
negative. The case $n,m>0$ corresponds to a subsequence of the form $\cdots
a_{1}vb^{\alpha}a_{1}v\cdots$ or $\cdots a_{1}b^{\alpha}v^{\ast}a_{1}v\cdots$.
In either case we use the fact that $E_{B_{0}}\left(  b^{\alpha}v\right)  =0$,
$E_{B_{0}}\left(  v^{\ast}b^{\alpha}\right)  =0$, for $\alpha\neq0$. The case
$n,m<0$ is similar.

Hence the property of the amalgamated free product applies, and $W$ is
non-scalar.\medskip

\emph{Step} IV. Any word (except the identity) in the Baumslag group
$\left\langle a,b\mid ab^{3}a^{-1}=b^{2}\right\rangle $, of total degree zero
in $a$, is equal to one of the words
\[
a^{\alpha_{1}}b^{\beta_{1}}a^{\alpha_{2}}b^{\beta_{2}}\cdots b^{\beta_{n}%
}a^{\alpha_{n+1}}%
\]
for which all the Assumptions A1--A3 on consecutive indices, described in Step
III, hold. Note that by Proposition \ref{proposition2.3}$\mathrm{(f)}$ we can
reduce the proof of the theorem to words of total degree $0$ in $a$.\medskip

To prove the claim of Step IV, the following two lemmas, dealing with easier
situations, will be used.

\begin{lemma}
\label{lemma2.3}Let $n\geq1$ and $k\geq2$. Then $a^{-n}b^{k}a^{n}$ is equal to
a product of the form
\[
b^{\theta_{0}}a^{-i_{1}}b^{\theta_{1}}a^{-i_{2}}\cdots a^{-i_{p}}b^{\theta
_{p}}a^{i_{p+1}}%
\]
for some strictly positive numbers $i_{1},\dots,i_{p+1}\in\{1,2,\dots,n\}$,
$i_{p+1}\leq n$, and $\theta_{0}$ is nonzero, $\left\vert \theta
_{0}\right\vert \geq3$, $\theta_{1},\dots,\theta_{p}=1$.

For example $a^{-n}b^{2}a^{n}=b^{3}a^{-1}ba^{-1}b\cdots ba^{-1}ba^{n-1}$,
where the product involves $n$ occurrences of the letter $b$ \textup{(}not
counting powers\/\textup{).}
\end{lemma}

\begin{proof}
We start with $k=2^{l}\cdot q$, $q$ odd. Then $a^{-n}b^{k}a^{n}=a^{-(n-l)}%
b^{3^{l}q}a^{n-l}$. We then split this product as:
\[
\left(  a^{-(n-l)}b^{3^{l}q-1}a^{n-l}\right)  \cdot\left(  a^{-(n-l)}%
ba^{(n-l)}\right)  .
\]
We repeat this procedure with
\[
3^{l}q-1=2^{l_{1}}q_{1},\text{\qquad}q_{1}\text{ odd }%
\]
and will obtain
\[
a^{-n}b^{k}a^{n}=a^{-(n-l-l_{1})}b^{3^{l_{1}}q_{1}}a^{(n-l-l_{1})}%
a^{-(n-l)}ba^{(n-l)}.
\]
By repeating this procedure and stopping when running out of powers of $a$, we
get the required result.
\end{proof}

Similarly one proves:

\begin{lemma}
\label{lemma2.4}Let $n\geq1$, $k\geq3$. Then $a^{n}b^{k}a^{-n}$ is equal to a
product of the form
\[
b^{\theta_{0}}a^{\varepsilon_{1}}b^{\theta_{1}}a^{\varepsilon_{2}}%
b^{\theta_{2}}\cdots a^{\varepsilon_{s}}b^{\theta_{s}}a^{-\varepsilon_{s+1}},
\]
where $0<\varepsilon_{s+1}\leq n$, $\varepsilon_{1},\varepsilon_{2}%
,\dots,\varepsilon_{s}>0$. Moreover, $\theta_{0}\geq2$ and $\theta_{s}%
\in\left\{  1,2\right\}  $ for $s\geq1$.
\end{lemma}

\begin{proof}
[Proof of the claim in Step \textup{IV}]We start with one arbitrary word
\[
W=b^{k_{0}}a^{n_{1}}b^{k_{1}}a^{n_{2}}b^{k_{2}}\cdots,\qquad n_{i}%
\neq0,\;k_{i}\neq0,\;i\geq1,
\]
where no obvious cancellations are possible. By moving from the left to the
right we look at the first change in sign in the sequence $n_{1},n_{2},\dots$.
Say this occurs when $i=i_{0}$. At that point, if $\left\vert k_{i_{0}%
}\right\vert \geq2$ when $n_{i_{0}}<0$, $n_{i_{0}+1}>0$, or if $\left\vert
k_{i_{0}}\right\vert \geq3$ when $n_{i_{0}}>0$, $n_{i_{0}+1}<0$, we apply one
of the two preceding lemmas to replace $a^{n_{i_{0}}}b^{k_{i_{0}}}%
a^{n_{i_{0}+1}}$ by one of the sequences described in the lemmas. More
precisely, if, for example, $n_{i_{0}}<0$, $n_{i_{0}+1}>0$, we apply Lemma
\ref{lemma2.3} for
\[
x=a^{-\min\left(  -n_{i_{0}},n_{i_{0}+1}\right)  }b^{k_{i_{0}}}a^{\min\left(
-n_{i_{0}},n_{i_{0}+1}\right)  }.
\]
By replacing $x$ in the word by the form given in Lemma \ref{lemma2.3}, the
structure of the word up to the next power of $a$ following $b^{k_{i_{0}}}$
would fulfill the requirements of the claim. The only case, when in doing this
replacement, a change of structure could occur in the structure of the word,
before $a^{n_{i_{0}}}$, is when $n_{i_{0}-1}>0$. But in this case $\left\vert
k_{i_{0}-1}\right\vert \leq2$ so $b^{k_{i_{0}-1}}$ won't cancel the $b^{3}$
appearing at the beginning of the word from Lemma \ref{lemma2.3}. Here we
reiterate the procedure. A similar argument works for $n_{i_{0}}>0$,
$n_{i_{0}+1}<0$.\renewcommand{\qed}{}
\end{proof}

By induction, this completes the proof of Step IV. By Steps III and IV we
conclude the proof of our theorem.
\end{proof}

\section{Extremal finite von Neumann algebras}

In this section we consider the structure of the set of moments of families of
projections in a finite von Neumann algebra. Note that by Kirchberg's
technique \cite{Ki}, for Connes's conjecture to be true, one should prove that
the closure of this set is independent of the finite von Neumann algebra for
which we consider the set of moments.

\begin{definition}
\label{definition3.1}Let $M$ be a finite separable von Neumann algebra and let
$\tau$ be a faithful, normalized trace on $M$. For any integer $n\geq1$, let
$K_{M}^{n}$ be the subset of $[0,1]^{\frac{n(n+1)}{2}}$ consisting of the
following ordered pairs:
\[
K_{M}^{n}=\{(\tau(e_{i}e_{j}))_{1\leq i\leq j\leq n}\mid(e_{1},e_{2}%
,\dots,e_{n})\in(\mathcal{P}(M))^{n}\}.
\]

\end{definition}

\begin{proposition}
\label{remark3.2}Let $M$ be a type $\mathrm{II}_{1}$ factor with trace $\tau$.
Then for all integers $n\geq1$,

\begin{enumerate}
\item $K_{M}^{n}$ is convex if $\mathcal{F}(M)=\mathbb{R}_{+}\backslash\{0\}$,

\item $K_{M}^{n}$ is closed under pointwise multiplication, if $M\cong
M\otimes M$,

\item $K_{M^{\omega}}^{n}$ is closed in the standard topology of
$[0,1]^{\frac{n(n+1)}{2}}$,

\item $K_{M}^{n}$ $\supseteq K_{R}^{n}$, where $R$ is the hyperfinite
$\mathrm{II}_{1}$ factor.
\end{enumerate}
\end{proposition}

The proof of this proposition is identical to the proof of the properties for
the set of moments associated with unitaries in a $\mathrm{II}_{1}$ factor.
Note that by Kirchberg's results \cite{Ki}, $\overline{K_{M}^{n}}%
=\overline{K_{R}^{n}}$ for all $n$, if and only if $M\subseteq R^{\omega}$.

It is very easy to describe the geometry of a diffuse abelian von Neumann
algebra. Indeed,

\begin{proposition}
\label{remark3.3}Let $Y_{n}\subseteq\lbrack0,1]^{\frac{n(n+1)}{2}}$ consist of
all $(\varepsilon_{ij})_{1\leq i\leq j\leq n}$ such that there are sets
$A_{1},\dots,A_{n}\subseteq X$, $X$ nonvoid, $A_{i}=\emptyset$ or $A_{i}=X$
such that $\varepsilon_{ij}=1$ if $A_{i}\cap A_{j}=X$ and $\varepsilon_{ij}=0$
if $A_{i}\cap A_{j}=\emptyset$. Then
\[
K_{L^{\infty}([0,1])}^{n}=\operatorname*{co}Y_{n}.
\]

\end{proposition}

In this section we analyze the structure of the closed convex subsets
$K_{M^{\omega}}^{n}\subseteq\lbrack0,1]^{\frac{n(n+1\}}{2}}$. To determine
completely this set it would be sufficient to know, for all choices of real
numbers $(a_{ij})_{1\leq i\leq j\leq n}$ of the value of
\[
\underset{1\leq i\leq j\leq n}{\max}\left\{  \sum a_{ij}\lambda_{ij}%
\mid(\lambda_{ij})\in K_{M}^{n}\right\}  .
\]
This is difficult to handle, but we are able to prove at least one geometrical
property related to this maximum value: a type of separation of variables at
maximum points in $K_{M}^{n}$.

The following lemma is an easy consequence of the fact that whenever a maximum
point is attained at $(e_{1}^{0},\dots,e_{n}^{0})$, then for any other
projection $e_{1}\leq e_{1}^{0}$ or $e_{1}\geq e_{1}^{0}$ we get a lower value.

\begin{lemma}
\label{lemma3.4}Fix
\[
(\tau(e_{i}^{0}e_{j}^{0}))_{1\leq i\leq j\leq n}%
\]
a maximum point for the fixed functional
\[
L(\lambda_{ij})=\sum_{1\leq i\leq j\leq n}a_{ij}\lambda_{ij}\text{ \qquad on
}K_{M}^{n}.
\]
For all $i=1,2,\dots,n$, let
\[
\Omega_{i}=\Omega_{i}(e_{1}^{0},\dots,e_{0}^{n})=\sum_{j\neq i}a_{ij}e_{j}%
^{0}+a_{ii}e_{i}^{0}.
\]
Then $e_{i}^{0}\Omega_{i}e_{i}^{0}\geq0$ and $(1-e_{i}^{0})\Omega_{i}%
(1-e_{i}^{0})\leq0$ for all $i=1,2,\dots,n$.
\end{lemma}

\begin{proof}
Fix $i$ in $\{1,2,\dots,n\}$ and let $e_{i}$ be any projection less than
$e_{i}^{0}$. The fact that
\[
\sum_{1\leq i\leq j\leq n}a_{ij}\tau\left(  e_{i}^{0}e_{j}^{0}\right)
\]
is a maximum value for $L$ on $K_{M}^{n}$, implies that
\[
\sum_{j\neq i}a_{ij}\tau\left(  (e_{i}^{0}-e_{i})e_{j}^{0}\right)  +a_{ii}%
\tau(e_{i}^{0}-e_{i})\geq0.
\]
Thus for any projection $e$ less than $e_{i}^{0}$ we have that%
\[
\tau\left(  e\left(  \sum_{j\neq i}a_{ij}e_{j}+a_{ii}\operatorname*{Id}%
\right)  \right)  \geq0.
\]
But this gives exactly that
\[
e_{i}^{0}\Omega_{i}e_{i}^{0}\geq0.
\]
Similarly for $1-e_{i}^{0}$.
\end{proof}

\begin{corollary}
\label{corollary3.5}If $\left(  \lambda_{ij}^{0}\right)  _{1\leq i\leq j\leq
n}$ in $K_{M}^{n}$ is a maximum point for
\[
\left(  \lambda_{ij}\right)  _{1\leq i\leq j\leq n}\in K_{M}^{n}%
\longrightarrow\sum_{1\leq i\leq j\leq n}a_{ij}\lambda_{ij},
\]
then for all $i=1,2,\dots,n$ we have that
\[
0\leq\sum_{j\neq i}a_{ij}\lambda_{ij}^{0}+a_{ii}\lambda_{ii}^{0}\leq
\sum_{j\neq i}a_{ij}\lambda_{jj}^{0}+a_{ii}.
\]

\end{corollary}

\begin{proof}
This follows by writing down explicitly that
\[
\tau(e_{i}^{0}\Omega_{i})\geq0\text{,\qquad}\tau\left(  \left(  1-e_{i}%
^{0}\right)  \Omega_{i}\right)  \leq0.
\]

\end{proof}

We will use a method similar to the method of Lagrange multipliers to
determine the finer structure of a set of projections $e_{i}^{0},\dots
,e_{n}^{0}$ at which a maximum point is attained.

To do this we need to show that the grassmanian manifold associated with a
type $\mathrm{II}_{1}$ factor is large enough.

\begin{lemma}
\label{lemma3.6}Let $M$ be a $\mathrm{II}_{1}$ factor, $e$ be a non-trivial
projection in $M$ and $\mathcal{T}_{e}$ be the linear space consisting of all
$Z$ in $M$, such that $Z=Z^{\ast}$ and $eZe=0$, $(1-e)Z(1-e)$. Let
$\overset{\circ}{\mathcal{T}}_{e}$ be the set of all $Z$ in $\mathcal{T}_{e}$
such that there exists a one-parameter family $e(t)$ of projections in $M$,
weakly differentiable at $0$, such that
\[
e(0)=e,\qquad\dot{e}(0)=Z.
\]
Then the space of $\overset{\circ}{\mathcal{T}}_{e}$ is weakly dense in
$\mathcal{T}_{e}$.
\end{lemma}

\begin{proof}
Assume first that $\tau(e)=\frac{1}{2}$, and let $v$ be any partial isometry
mapping $e$ onto $1-e$. We will show that $Z=v+v^{\ast}$ belongs to
$\overset{\circ}{\mathcal{T}}_{e}$.

Indeed $\{e,v\}^{\prime\prime}$ can be identified with $M_{2}(\mathbb{C})$ in
such a way that
\[
v+v^{\ast}=\left(
\begin{array}
[c]{ll}%
0 & 1\\
1 & 0
\end{array}
\right)  ,\text{\qquad}e=\left(
\begin{array}
[c]{ll}%
0 & 0\\
0 & 1
\end{array}
\right)  .
\]
But then we take
\[
e(\theta)=\left(
\begin{array}
[c]{cc}%
\sin^{2}\theta & \sin\theta\cos\theta\\
\sin\theta\cos\theta & \cos^{2}\theta
\end{array}
\right)  .
\]
If the trace of $e$ is different from $\frac{1}{2}$, we may then assume that
$\tau(e)<\frac{1}{2}$. By the above argument, any partial isometry $v$ mapping
$e$ into a projection under $1-e$, determines an element $v+v^{\ast}$ in
$\overset{\circ}{\mathcal{T}}_{e}$.

Thus $\overset{\circ}{\mathcal{T}}_{e}$ contains $v+v^{\ast}$ for any partial
isometry $v$, such that $v^{\ast}v=e$, $vv^{\ast}\leq1-e$. Let $u$ be any
unitary in $eMe$ and let $w$ be any unitary in $\left(  1-e\right)  M\left(
1-e\right)  $. The same argument shows that $wvu+u^{\ast}v^{\ast}w^{\ast}$
belongs to $\overset{\circ}{\mathcal{T}}_{e}$. Since any element in $eMe$ and
$\left(  1-e\right)  M\left(  1-e\right)  $ is a linear combination of
unitaries, this shows that $yvx+x^{\ast}v^{\ast}y^{\ast}$ is always in the
linear span of $\overset{\circ}{\mathcal{T}}_{e}$ for all $x$ in $eMe$ and $y$
in $\left(  1-e\right)  M\left(  1-e\right)  $. This set is obviously weakly
dense in $\mathcal{T}_{e}$.
\end{proof}

\begin{corollary}
\label{corollary3.7}Fix $n\geq1$ and real numbers $(a_{ij})_{1\leq i\leq j\leq
n}$. Let
\[
(e_{1}^{0},e_{2}^{0},\dots,e_{n}^{0})
\]
be a family of projections in
$\mathcal{P}(M)$ such that the maximum of
\[
L\left(  \left(  \lambda_{ij}\right)  _{1\leq i\leq j\leq n}\right)
=\sum_{1\leq i\leq j\leq n}a_{ij}\lambda_{ij}%
\]
for $\left(  \lambda_{ij}\right)  $ in $K_{M}^{n}$ is attained at $(\tau
(e_{i}^{\circ}e_{j}^{\circ}))$. Let
\[
\Omega_{i}^{0}=\Omega_{i}^{0}(e_{i}^{0},\dots,e_{n}^{0},a_{ij})=a_{ii}%
\operatorname*{Id}+\sum_{j\neq i}a_{ij}e_{j}^{0}.
\]
Then $[e_{i}^{0},\Omega_{i}^{0}]=0$. By using Lemma $\ref{lemma3.4}$ it
follows that
\begin{align*}
e_{i}^{0}  &  \geq\operatorname*{supp}(\Omega_{i}^{0})_{+},\\
1-e_{i}^{0}  &  \geq\operatorname*{supp}\left(  \Omega_{i}^{0}\right)  _{-}.
\end{align*}

\end{corollary}

\begin{proof}
Indeed if
\[
(\tau(e_{i}^{0}e_{j}^{0}))_{1\leq i\leq j\leq n}%
\]
is such a maximum point for the functional $L$ on $K_{M}^{n}$, then for all
$Z_{i}$ in $\overset{\circ}{\mathcal{T}}_{e_{i}^{0}}$ we have that
\[
\tau(\Omega_{i}^{0}Z_{i})=0.
\]
But then this will give that
\[
\tau(\Omega_{i}^{0}Z_{i})=0
\]
for all $Z_{i}=(1-e_{i}^{0})Y_{1}e_{i}^{0}+e_{i}^{0}Y_{1}(1-e_{i}^{0})$,
$Y_{1}\in M_{sa}$. Thus for all $Y=Y^{\ast}$ in $M$ we have
\[
\tau\left(  \Omega_{i}^{0}e_{i}^{0}Y(1-e_{i}^{0})+\Omega_{i}^{0}(1-e_{i}%
^{0})Ye_{i}^{0}\right)  =0
\]
and hence%
\[
\tau\left(  \lbrack(1-e_{i}^{0})\Omega_{i}^{0}e_{i}^{0}+e_{i}^{0}\Omega
_{i}^{0}(1-e_{i}^{0})]Y\right)  =0
\]
for all $Y$ selfadjoint in $M$. Since $\Omega_{i}^{0}$ is also selfadjoint, it
follows that
\[
(1-e_{i}^{0})\Omega_{i}^{0}e_{i}^{0}+(1-e_{i}^{0})\Omega_{i}^{0}e_{i}^{0}=0
\]
and hence that $\Omega_{i}^{0}$ commutes with $e_{i}^{0}$.
\end{proof}

\begin{remark}
The above proposition suggests that for Connes's embedding problem, it is
sufficient to consider finite von Neumann algebras \textup{(}which we call
\emph{extremal finite von Neumann algebras\/}\textup{)} that are generated by
families of projections $e_{1},e_{2},\dots,e_{n}$ such that there exists a
matrix of real numbers $(a_{ij})_{1\leq i\leq j\leq n}$ with the following property.

For each $i$, let
\[
\Omega_{i}=\sum_{j\neq i}a_{ij}e_{j}+a_{ii}\operatorname*{Id}.
\]
Let $s_{+}^{i}$ be the init of the positive part of $(\Omega_{i})_{+}$ and
$s_{-}^{i}$ be the projection onto the init space of $(\Omega_{i})_{-}$.

Then $1-s_{-}^{i}\geq e_{i}\geq s_{+}^{i}$; in particular, $e_{i}$ commutes
with $\Omega_{i}$, for all $i$.
\end{remark}

\ifx\undefined\bysame\fi

\end{document}